\documentstyle [12pt,epsf, graphicx, amssymb]{article}

\begin{document}
\def\l{\lambda}
\def\m{\mu}
\def\a{\alpha}
\def\b{\beta}
\def\g{\gamma}
\def\d{\delta}
\def\e{\epsilon}
\def\o{\omega}
\def\O{\Omega}
\def\v{\varphi}
\def\t{\theta}
\def\r{\rho}
\def\bs{$\blacksquare$}
\def\bp{\begin{proposition}}
\def\ep{\end{proposition}}
\def\bt{\begin{th}}
\def\et{\end{th}}
\def\be{\begin{equation}}
\def\ee{\end{equation}}
\def\bl{\begin{lemma}}
\def\el{\end{lemma}}
\def\bc{\begin{corollary}}
\def\ec{\end{corollary}}
\def\pr{\noindent{\bf Proof: }}
\def\note{\noindent{\bf Note. }}
\def\bd{\begin{definition}}
\def\ed{\end{definition}}
\def\C{{\mathbb C}}
\def\P{{\mathbb P}}
\def\Z{{\mathbb Z}}
\def\d{{\rm d}}
\def\deg{{\rm deg\,}}
\def\deg{{\rm deg\,}}
\def\arg{{\rm arg\,}}
\def\min{{\rm min\,}}
\def\max{{\rm max\,}}

\newcommand{\norm}[1]{\left\Vert#1\right\Vert}
\newcommand{\abs}[1]{\left\vert#1\right\vert}

\newcommand{\set}[1]{\left\{#1\right\}}
\newcommand{\setb}[2]{ \left\{#1 \ \Big| \ #2 \right\} }

\newcommand{\IP}[1]{\left<#1\right>}
\newcommand{\Bracket}[1]{\left[#1\right]}
\newcommand{\Soger}[1]{\left(#1\right)}

\newcommand{\Integer}{\mathbb{Z}}
\newcommand{\Rational}{\mathbb{Q}}
\newcommand{\Real}{\mathbb{R}}
\newcommand{\Complex}{\mathbb{C}}

\newcommand{\eps}{\varepsilon}
\newcommand{\To}{\longrightarrow}
\newcommand{\varchi}{\raisebox{2pt}{$\chi$}}

\newcommand{\E}{\mathbf{E}}
\newcommand{\Var}{\mathrm{var}}

\def\squareforqed{\hbox{\rlap{$\sqcap$}$\sqcup$}}
\def\qed{\ifmmode\squareforqed\else{\unskip\nobreak\hfil
\penalty50\hskip1em\null\nobreak\hfil\squareforqed
\parfillskip=0pt\finalhyphendemerits=0\endgraf}\fi}

\renewcommand{\th}{^{\mathrm{th}}}
\newcommand{\Dif}{\mathrm{D_{if}}}
\newcommand{\Difp}{\mathrm{D^p_{if}}}
\newcommand{\GHF}{\mathrm{G_{HF}}}
\newcommand{\GHFP}{\mathrm{G^p_{HF}}}
\newcommand{\f}{\mathrm{f}}
\newcommand{\fgh}{\mathrm{f_{gh}}}
\newcommand{\T}{\mathrm{T}}
\newcommand{\K}{^\mathrm{K}}
\newcommand{\PghK}{\mathrm{P^K_{f_{gh}}}}
\newcommand{\Dig}{\mathrm{D_{ig}}}
\newcommand{\for}{\mathrm{for}}
\newcommand{\End}{\mathrm{end}}

\newtheorem{th}{Theorem}[section]
\newtheorem{lemma}{Lemma}[section]
\newtheorem{definition}{Definition}[section]
\newtheorem{corollary}{Corollary}
\newtheorem{proposition}{Proposition}[section]

\begin{titlepage}

\begin{center}

\topskip 5mm

{\LARGE{\bf {Generalized Remez Inequality for $(s,p)$-Valent Functions}}} \vskip 8mm

{\large {\bf Y. Yomdin}}

\vspace{6 mm}

{Department of Mathematics, The Weizmann Institute of Science,
Rehovot 76100, Israel. e-mail: yosef.yomdin@weizmann.ac.il}

\end{center}

\vspace{2 mm}
\begin{center}

{ \bf Abstract}
\end{center}

{\small The classical Remez inequality bounds the maximum of the
absolute value of a polynomial $P(x)$ of degree $d$ on $[-1,1]$
through the maximum of its absolute value on any subset $Z$ of
positive measure in $[-1,1]$. It was shown in \cite{Yom3} that
the Lebesgue measure in the Remez inequality can be replaced by a
certain geometric invariant $\omega_d(Z)$ which can be effectively
estimated in terms of the metric entropy of $Z$ and which may be
nonzero for discrete and even finite sets $Z$.

In the present paper we first obtain an essentially sharp Remez-type
inequality in the spirit of \cite{Yom3} for complex polynomials of
one variable, introducing metric invariants
$c_d(Z)$ and $\o_{cd}(Z)$ for an arbitrary subset
$Z\subset D_1\subset {\mathbb C}$. These invariants translate into the
the metric language the classical Cartan lemma (see \cite{Gor} and
references therein).

Next we introduce $(s,p)$-valent functions, which provide a natural
generalization of $p$-valent ones (see \cite{Hay} and references
therein). We prove a ``distortion theorem" for such functions, comparing
them with polynomials sharing their zeroes. On this base we extend to
$(s,p)$-valent functions our polynomial Remez-type inequality.

As the main example we consider restrictions $g$ of polynomials of a
growing degree to a fixed algebraic curve, for which
we obtain an essentially sharp ``local" Remez-type inequality,
stressing the role of the geometry of singularities of $g$.

Finally, we obtain for such functions $g$ a ``global" Remez-type
inequality which is valid for all the branches of $g$ and involves
both the geometry of singularities of $g$ and its monodromy.}

\vspace{2 mm}
\begin{center}
------------------------------------------------
\vspace{2 mm}
\end{center}
This research was supported by the ISF, Grant No. 304/09, by the Minerva
Foundation, and by the European Social Fund.

\end{titlepage}

\newpage

\section{Introduction}
\setcounter{equation}{0}

The classical Remez inequality bounds the maximum of the absolute
value of a polynomial $P(x)$ of degree $d$ on $[-1,1]$ through the
maximum of its absolute value on any subset $Z$ of positive
measure in $[-1,1]$. More accurately, we have (\cite{Rem}):

\medskip

\noindent {\it Let $P(x)$ be a real polynomial of degree $d$. Then for
any measurable $Z\subset [-1,1]$ \be \max_{[-1,1]} \vert P(x)
\vert \leq T_d({{4-\mu}\over {\mu}})\max_Z \vert P(x) \vert ,\ee
where $\mu=\mu_1(Z)$ is the Lebesgue measure of $Z$ and
$T_d(x)=cos(d \ arccos(x))$ is the $d$-th Chebyshev polynomial}.

\smallskip

However, the inequality of the form (1.1) may be true also for
some sets $Z$ of measure zero and even for certain finite sets
$Z$. In particular, in \cite{Fav} certain ``fractal" subsets
$Z\subset [-1,1]$ are studied. An invariant $\phi_Z(d)$ is defined
and estimated in some examples, which is the best constant in the
Remez inequality for the couple $(Z\subset [-1,1])$.

\smallskip

In \cite{Yom3} it was shown  that the Lebesgue measure in the Remez
inequality (in one and several variables) can be replaced by a certain
geometric invariant $\omega_d(Z)$ which can be effectively estimated
in terms of the metric entropy of $Z$ and which may be nonzero for
discrete and even finite sets $Z$. In one-dimensional case the definition
of $\omega_d(Z)$ is very simple: $\omega_d(Z)= sup_{\e>0} \ \e(M(Z,\e)-d),$
where $M(Z,\e)$ is a minimal number of $\e$-intervals covering $Z$.
So we have a generalized Remez inequality of the following form:

\bt (\cite{Yom3}) Let $P(x)$ be a real polynomial of degree $d$. Then for
any $Z\subset [-1,1]$ \be \max_{[-1,1]} \vert P(x)
\vert \leq T_d({{4-\o}\over {\o}})\max_Z \vert P(x) \vert ,\ee
where $\o=\o_d(Z)$ and $T_d(x)=cos(d \ arccos(x))$ is the $d$-th
Chebyshev polynomial. \et Since for a measurable $Z$ we have clearly
$\o_d(Z) \geq \mu_1(Z)$ (let $\e$ tend to zero) Theorem 1.1
is a true generalization of the classical Remez bound.

\smallskip

The aim of the present paper is to extend the result of Theorem 1.1
to complex analytic functions of one variable. A lot of
inequalities in this direction (in one and several variables) are known
(see, for example, \cite{Bru2,Bru.Bru1,Bru.Bru2,Zer1,Zer2}
and references therein). However, to our best knowledge, there are still
some important open questions concerning Remez-type inequalities for
complex analytic functions, which we answer in this paper:

\smallskip

1. It would be desirable to have a simple geometric invariant $c_d(Z)$
computable for every subset $Z$ of a complex unit disk $D_1$, and
providing a nontrivial information also for finite and discrete sets $Z$,
and an inequality of the form \be \max_{D_1} \vert P(x) \vert \leq
\Psi_d(c_d(Z))\max_{Z} \vert P(x) \vert \ee for each complex polynomial
$P(x)$ of degree $d$. {\it Below we provide both the invariant $c_d$ and
the inequality of the form (1.3) based on the classical Cartan (or Cartan -
Boutroux) lemma.}

\smallskip

2. As we want to replace polynomials in (1.3) by more general analytic
functions $f$, it is not clear from the results available
today, what information on $f$ is truly relevant. Typically, Remez-type
inequalities are studied for some special finite-dimensional families, 
like exponential polynomials or quasi-polynomials (\cite{Tur,Naz,Bru1}.

In this paper we introduce $(s,p)$-valent functions in a domain $\O$,
which provide a natural generalization of classically known $p$-valent
functions. The last are those functions $f$ for which the equation
$f(x)=c$ has at most $p$ solutions in $\O$ for any $c$ (see \cite{Hay}
and references therein). For an $(s,p)$-valent function $f$ the equation
$f(x)=P(x)$ has at most $p$ solutions in $\O$ for any polynomial $P(x)$
of degree $s$. $(0,p)$-valent functions are the same as $p$-valent.

There are many natural examples of $(s,p)$-valent functions: algebraic
functions, and, in particular, restrictions of polynomials of a growing
degree to a fixed algebraic curve, which we consider below, belong to
this class. The same is true for solutions of algebraic differential
equations. However, as we show in \cite{Roy.Yom2}, this notion is
applicable to any analytic function, under an appropriate choice of
the domain $\O$ and the parameters $(s,p)$, and it may provide a useful
information in very general situations.

We start here an investigation of $(s,p)$-valent functions, following
the classically known patterns for the $p$-valent ones. In particular,
{\it we prove (following \cite{Hay}) a pretty accurate ``distortion
theorem" for such functions, which compares them with the polynomials
sharing their zeroes.

On this base we extend to $(s,p)$-valent functions our polynomial
Remez-type inequality (1.3)}.

\smallskip

3. As the main example we consider restrictions $g$ of polynomials of a
growing degree to a fixed algebraic curve. It turns out that the only
information we need in order to extend to this case the inequality of
the form (1.3) are the degrees of the curve and of the restricted
polynomial, and the distance to the nearest singularity.

Finally, we obtain for such functions $g$ a {\it ``global" Remez-type
inequality which is valid for all the branches of $g$. It involves
the geometry of singularities of $g$, and its monodromy.}

\medskip

The paper is organized as follows: in Section 2 we define the
invariant $c_d$ and prove a ``discrete Remez inequality" of the form (1.3)
for complex polynomials. In Section 3 we introduce $(s,p)$-valent functions
and prove a ``distortion theorem". In Section 4 we obtain Remez inequality
of the form (1.3) for $(s,p)$-valent functions. Finally, in Section 5 we
obtain ``local" and ``global" Remez inequalities for restrictions of
polynomials to an algebraic curve.

\medskip

The author would like to express his gratitude to M. Briskin and
A. Brudnyi for useful discussions, and to the Center for Advanced
Studies, Warsaw University of Technology, for its support and kind
hospitality.

\section{Remez inequality for complex polynomials}
\setcounter{equation}{0}

For polynomials in one complex variable a result similar to the
Remez inequality is provided by the classical Cartan (or Cartan -
Boutroux) lemma (see, for example, \cite{Gor} and references therein):

\smallskip

\noindent {\it Let $P(z)$ be a monic polynomial of one complex
variable of degree $d$. For any given $\e>0$ consider
$V_{\e^d}(P)=\{z\in {\mathbb C}, \ \vert P(z) \vert \leq \e^d\}.$
Then $V_{\e^d}(P)$ can be covered by $d$ complex discs $D_j$ with
radii $r_j, \ j=1,\dots,d,$ such that $\sum_{j=1}^d r_j \leq 2e
\e.$}

\smallskip

In \cite{Bru,Bru.Bru2,Zer2,Zer1} some generalizations of the
Cartan - Boutroux lemma to plurisubharmonic functions have been
obtained, which lead, in particular, to the bounds on the size
of sub-level sets. In these lines in \cite{Bru} some bounds for
the covering number of sublevel sets of complex analytic functions
have been obtained, similar to the results of \cite{Yom3} in the
real case.

\medskip

In the present paper we would like to derive from the Cartan lemma
both the definition of the invariant $c_d$ and the corresponding Remez
inequality.
\bd For any $Z\subset D_1$ let $c_d(Z)$ be the minimal sum of the radii
of d disks covering $Z$.\ed Now we can state and proof our generalized
Remez inequality for complex polynomials:

\bt For $Z\subset D_1$ assume that $c_d(Z)=c>0$. Then each complex
polynomial $P(x)$ of degree $d$ satisfies:
\be \max_{D_1} \vert P(x) \vert \leq ({6e \over c})^d \ \max_{Z}
\vert P(x) \vert.\ee \et\pr We obtain the proof in two steps. First we
prove a lemma which is just a reformulation of the Cartan's one:
\bl Assume that $\vert P(x) \vert \leq 1$
on $Z$. Then the leading coefficient of $P$ does not exceed in absolute
value $({2e \over c})^d$.\el\pr For each $\rho \geq 0$ put $V_\rho(P) =
\{z\in {\mathbb C}, \ \vert P(z) \vert \leq \rho\}.$ Since by assumptions
$\vert P(x) \vert \leq 1$ on $Z$ we have $Z\subset V_{1}(P).$ By the
definition of $c_d(Z)$ we conclude that for every covering of
$V_{1}(P)$ by $d$ disks $D_1,\dots,D_d$ of the radii $r_1,\dots,r_d$
(which is also a covering of $Z$) we have
$\sum_{i=1}^d r_i \geq c_d(Z)$. Denoting the absolute value of the leading
coefficient of $P(x)$ by $A$ we have by the Cartan lemma that for a
certain covering as above
$c=c_d(Z)\leq \sum_{i=1}^d r_i \leq 2e({1\over A})^{1\over d}.$
We conclude that $A\leq ({2e \over c})^d$. $\square$

\smallskip

Now we write $P(x)=A\prod_{j=1}^d (x-x_j), \ x_j \in {\mathbb C},$ and
consider separately two cases:

\smallskip

1. $\vert x_j \vert \leq 2, \ j=1,\dots,d.$ In this case
$\max_{D_1} \vert P(x) \vert \leq 3^dA\leq 3^d ({2e \over c})^d$. So in
this case the proof of Theorem 2.1 is completed.

\smallskip

2. $\vert x_j \vert \leq 2, \ j=1,\dots,d_1<d,$ while
$\vert x_j \vert > 2, \ j=d_1+1,\dots,d.$ In this case we write
$P_1(x)=A\prod_{j=1}^{d_1} (x-x_j), \ P_2(x)=\prod_{j=d_1+1}^d (x-x_j),$
and notice that for any two points $v_1,v_2\in D_1$ the ratio of absolute
values of $P_2(v_1)$ and $P_2(v_1)$ satisfies
$\vert {{P_2(v_1)}\over {P_2(v_2)}}\vert < 3.$ Consequently we get
\be{{\max_{D_1} \vert P(x) \vert}\over {\max_{Z} \vert P(x) \vert}}<
3^{d-d_1}{{\max_{D_1} \vert P_1(x)\vert}\over{\max_{Z} \vert P_1(x)\vert}}.\ee
All the roots of $P_1$ are bounded in absolute value by $2$, so by already
proved case of Theorem 2.1 we have
${{\max_{D_1} \vert P_1(x)\vert}\over{\max_{Z} \vert P_1(x)\vert}}\leq
3^{d_1}({2e \over c})^d$. Combining this inequality with (2.2) we get
once more ${{\max_{D_1} \vert P(x)\vert}\over{\max_{Z} \vert P(x)\vert}}
\leq 3^d({2e \over c})^d$. This completes the proof of Theorem 2.1.
$\square$

\smallskip

Clearly, the invariant $c_d(Z)$ may be positive for finite $Z$. In fact,
we have: \bp $c_d(Z) > 0$ if and only if Z contains more than $d$ points.\ep\pr
Indeed, any $d$ points can be covered by d disks with an arbitrarily small
sum of the radii. But the sum of radii of any $d$ disks covering at least
$d+1$ different points is greater than or equal to the one half of a minimal
distance between these points. $\square$

\smallskip

We now compare the invariant $c_d(Z)$ with some other
metric invariants which may be sometimes easier to compute. In particular,
we can easily produce a simple lower bound $c_d(Z)$  through the measure of $Z$:

\bp For each measurable $Z\subset D_1$ we have $c_d(Z)\geq
({{\mu_2(Z)}\over {\pi}})^{1/2}.$\ep\pr For the covering of $Z$ by
$d$ disks $D_1,\dots,D_d$ of the radii $r_1,\dots,r_d$ we have
$\sum_{i=0}^d r_i \geq (\sum_{i=0}^d r^2_i)^{1/2}.$ $\square$

\smallskip

However, in order to deal with discrete or finite subsets $Z\subset D_1$
we have to compare $c_d(Z)$ with the covering number $M(\e,Z)$ (which is,
by definition, the minimal number of $\e$-disks covering $Z$).

\bd For each $Z\subset D_1$ define $\o_{cd}(Z)$ as \ \
$\sup_\e \ \e[M(\e,Z)-d]^{1\over 2}.$ Define $\rho_d(Z)$ as \ \ $d\e_0$ where
$\e_0$ is the minimal $\e$ for which there is a covering of $Z$ with $d$
$\e$-disks (or $\e_0=M(.,Z)^{-1}(d)$). \ed {\bf Remark} The letter ``c"
in the notation of $c_d(Z)$ is for ``Cartan", and in the notation
of $\omega_{cd}(Z)$ it is for ``complex". As it was mentioned above, a very
similar invariant $\o_d(Z)=\sup_\e \ \e[M(\e,Z)-d]$ was introduced and used
in \cite{Yom3} in the real case. We compare $\omega_{cd}$ and $\o_d$ below.

\bt For each $Z\subset D_1$ we have
${1\over 4}\omega_{cd}(z)\leq c_d(Z)\leq \rho_d(Z)$.\et
\pr To prove the upper bound for $c_d(Z)$ we notice that it is the infinum
of the sum of the radii in all the coverings of $Z$ with $d$ disks, while
$\rho_d(Z)$ is such a sum for one specific covering.

To prove the lower bound, let us fix a covering of $Z$ by $d$ disks $D_i$
of the radii $r_i$
with $c_d(Z) = \sum_{i=0}^d r_i.$ Now for a given $\e > 0$ we conclude
that $M(\e,Z)$ is at most $d + ({{16}\over {\e^2}})\sum_{i=0}^d r^2_i.$
Indeed, for each disk $D_j$ with $r_j \geq \e$ we need at most
${{16r^2_j}\over {\e^2}}$ $\e$-disks to cover it. For each disk $D_i$
with $r_i \leq \e$ we need exactly one $\e$-disk to cover it, and the
number of such $D_i$ does not exceed $d$.
So we get $c_d(Z) \geq (\sum_{i=0}^d r^2_i)^{1\over 2} \geq {\e \over 4}
[M(\e,Z)-d]^{1\over 2}$. Taking supremum with respect to $\e>0$ we get
$c_c(Z)\geq {1\over 4}\omega_{cd}(Z).$ $\square$

\smallskip

Since $M(\e,Z)$ is always an integer, we have $\o_d(Z)\geq \o_{cd}(Z)$. For
$Z\subset D_1$ of positive plane measure $\o_d(Z)=\infty$ while $\o_{cd}(Z)$
remains bounded (in particular, by $\rho_d(Z)$). Some examples of computing
(or bounding) $\o_d(Z)$ for ``fractal" sets $Z$ can be found in \cite{Yom3}.
Computations for $\o_{cd}(Z)$ are essentially the same. In particular,
in an example given in \cite{Yom3} in connection to \cite{Fav} we
have: for
$Z=Z_r= \{1,{1\over {2^r}},{1\over {3^r}},\dots,{1\over {k^r}},\dots\}$
$$\o_d(Z_r) \asymp {{r^r}\over {(r+1)^{r+1}}} \ {1\over {d^r}}, \ \ \
\o_{cd}(Z_r) \asymp {{(2r+1)^r}\over {(2r+2)^{r+1}}} \ {1\over {d^{r+1/2}}} \ .$$
The asymptotic behavior here is for $d\rightarrow \infty$, as in \cite{Fav}.

\section{$(s,p)$-Valent Functions: Distortion Theorem}
\setcounter{equation}{0}

In this section we introduce a notion of an ``$(s,p)$-valent function"
which is a generalization of the classical notion of a $p$-valent function. 
The main reason is that for $(s,p)$-valent functions we can prove a kind
of ``distortion theorem" which shows that the behavior of $f$ is controlled
by the behavior of the polynomial with the same zeroes as $f$. Functions in
essentially all the classes traditionally studied in relation to
Bernstein-Markov-Remez type inequalities are $(s,p)$-valent, for any $s$ and
some $p$ depending on $s$, and we believe that there are good reasons to
study this property in general (see \cite{Roy.Yom2}.

\bd A function $f$ regular in a domain $\O\subset {\mathbb C}$ is called
$(s,p)$-valent in $\O$ if for any polynomial $P(x)$ of degree at most
$s$ the number of solutions of the equation $f(x)=P(x)$ in $\O$ does not
exceed $p$.\ed For $s=0$ we obtain the usual $p$-valent functions. Easy
examples show that an $(s,p)$-valent function may be not $(s+1,p)$-valent:

\bl The function $f(x)=x^p+ x^N, \ N\geq 10p+1,$ is $(s,p)$-valent in the
disk $D_{1\over 3}$ for any $s=0,\dots,p-1$, but only $(p,N)$-valent there.
\el\pr Taking $P(x)=x^p+c$ we see that the equation $f(x)=P(x)$ takes the
form $x^N = c$ so for small $c$ it has exactly $N$ solutions in
$D_{1\over 3}$. To show that $f$ is $(s,p)$-valent in the disk
$D_{1\over 3}$ for any $s=0,\dots,p-1$ fix a polynomial $P(x)$ of degree
$s\leq p-1$. Then the equation $f(x)=P(x)$ takes the form $-P(x)+x^p+x^N=0.$
Applying to the polynomial $Q(x)=-P(x)+x^p$ of degree $p$ (and with the leading
coefficient $1$) \ Lemma 3.3 of \cite{Yom2} we find a circle
$S_\rho = \{\vert x \vert = \rho\}$ with ${1\over 3}\leq \rho \leq {1\over 2}$
such that $\vert Q(x) \vert \geq ({1\over 2})^{10p}$ on $S_\rho$. On the other
hand $x^N \leq ({1\over 2})^{10p+1} < ({1\over 2})^{10p}$ on $S_\rho$.
Therefore by the Rouchet principle the number of zeroes of $Q(x)+x^N$ in
the disk $D_\rho$ is the same as for $Q(x)$, which is at most $p$. $\square$

\smallskip

A detailed study of the $(s,p)$-valent functions is presented in \cite{Roy.Yom2}.
In particular, we show there that solutions of algebraic differential
equations possess this property for any $s$ and some $p$ depending
on $s$. We also study in \cite{Roy.Yom2} wider classes of functions
whose Taylor coefficients obey ``algebraic" recurrence relations of a certain
form studied in \cite{Bri.Yom}, and show  that such functions are
$(s,p)$-valent for a proper choices of $s$ and $p$. Moreover, as we show in 
\cite{Roy.Yom2}, this notion is applicable to {\it any analytic function $f$}, 
under an appropriate choice of the domain $\O$ and the parameters $(s,p)$, with 
respect to the initial Taylor coefficients of $f$, and it may provide a useful
information in very general situations.

\medskip

Now we are ready to state and prove the ``distortion theorem" for
$(s,p)$-valent functions:

\bt Let $f(x)$ be a regular function in $D_1$ having there exactly $s$ zeroes
and being $(s,p)$-valent in $D_1$. Let $x_1,\dots, x_s$
be all the zeroes of $f$ in $D_1$, multiple zeroes counted according to their
multiplicity. Define $P(x)=a\prod_{j=1}^s(x-x_j)$, where the coefficient $a$ is
chosen in such a way that the constant term in the Taylor series for
$g(x)={{f(x)}\over {P(x)}}$ is equal to $1$. Then taking $\rho = \vert x \vert$
we have for any $x\in D_1$
\be ({{1-\rho}\over {1+\rho}})^{2p}\leq \vert g(x)
\vert \leq ({{1+\rho}\over {1-\rho}})^{2p}.\ee\et \pr The function
$g(x)={{f(x)}\over {P(x)}}$ is regular in $D_1$ and does not vanish there.
Moreover, $g$ is $p$-valent in $D_1$. Indeed, the equation $g(x)=c$ is equivalent
to $f(x)=cP(x)$ so it has at most $p$ solutions by definition of $(s,p)$-valent
functions. Now we apply Theorem 5.1 of \cite{Hay}: for a function $g$ as above
$({{1-\rho}\over {1+\rho}})^{2p}\leq \vert g(x)
\vert \leq ({{1+\rho}\over {1-\rho}})^{2p}.$ $\square$

\smallskip

It is not clear whether the requirement for $f$ to be $(s,p)$-valent is really
necessary in this theorem. In particular, in a special case where all the
zeroes of $f$ coincide, the following result is true:

\bt (Theorem 5.3, \cite{Hay}) Suppose that $f(x)=x^p+a_{p+1}x^{p+1}+\dots$
is $c.m. p$-valent in $D_1$. Then ${{\rho^p}\over ({1+\rho})^{2p}}\leq
\vert f(x) \vert \leq {{\rho}^p\over ({1-\rho})^{2p}}.$\et This theorem uses
a notion of ``circumferentially mean $p$-valent" ($c.m. p$-valent) functions,
which is weaker than $p$-valency. However, the
proof in \cite{Hay} is bases on a specific fact that under the assumptions
of the theorem $[f(x)]^{1\over p}$ turns out to be $m.p. 1$-valent in $D_1$.
It is not clear how to generalize this to zeroes of $f$ in a general position.
On the other hand, a ratio $g(x)={{f(x)}\over {P(x)}}$ certainly may be
not $p$-valent for $f$ being just $p$-valent, but not $(s,p)$-valent.
Indeed, take $f(x)=x^p+x^N$ as in Lemma 3.1 above. By this lemma $f$ is
$p$-valent in $D_{1\over 3}$ and it has a root of multiplicity $p$ at zero.
So $g(x)={{f(x)}\over {x^p}}=1+x^{N-p}$ and the equation $g(x)=c$ has
$N-p$ solutions in $D_{1\over 3}$ for $c$ sufficiently close to $1$. So
$g$ is not $p$-valent there.

\section{Remez Inequality for $(s,p)$-valent functions}
\setcounter{equation}{0}

The distortion theorem proved in the previous section allows us to
easily extend the Remez inequality of Section 2 from polynomials to
$(s,p)$-valent functions, just comparing them with polynomials having
the same zeroes.

\bt Let $f(x)$ be a regular function in $D_1$ having there exactly $s$ zeroes,
and being $(s,p)$-valent in $D_1$. Assume that $Z\subset D_R, \ R < 1,$ and
$c_d(Z)=c>0$. Then on each $D_{R'}$ with $0\leq R'<1$ $f$ satisfies
\be \max_{D_{R'}} \vert f(x) \vert \leq \sigma(R,R')^p({6e \over c})^s
\max_{Z}\vert f(x) \vert,\ee where
$\sigma(R,R')=[{{(1+R)(1+R')}\over {(1-R)(1-R')}}]^2.$\et\pr Assume that
$\vert f(x) \vert$ is bounded by $1$ on $Z$. Let $x_1,\dots, x_s$
be zeroes of $f$ in $D_1$. Consider, as in Theorem 3.1, the
polynomial $P(x)=a\prod_{j=1}^l(x-x_j)$, where the coefficient $a$ is
chosen in such a way that the constant term in the Taylor series for
$g(x)={{f(x)}\over {P(x)}}$ is equal to $1$. Then by Theorem 3.1 for
$g={f\over P}$ we have \be ({{1-\vert x \vert}\over {1+\vert x \vert}})^{2p}
\leq \vert g(x) \vert \leq ({{1+\vert x \vert}\over {1-\vert x \vert}})^{2p}.\ee
We conclude that
$P(x)\leq ({{1+R}\over {1-R}})^{2p}$ on $Z$. Hence by the polynomial
Remez inequality provided by Theorem 2.1 we obtain
$\vert P(x) \vert \leq ({6e \over c})^s({{1+R}\over {1-R}})^{2p}$ on $D_1$.
Finally, we apply once more the bound of Theorem 3.1 to conclude that
$\vert f(x) \vert \leq ({6e \over c})^s({{1+R}\over {1-R}})^{2p}
({{1+R'}\over {1-R'}})^{2p}$ on $D_{R'}$. $\square$

\section{Remez Inequality for polynomials on algebraic curves}
\setcounter{equation}{0}

Let $y=h(x)$ be an algebraic function, satisfying an equation
\be Q(x,h(x))=\sum_{i=0}^{d} s_i(x)h^i(x)=0,\ee with $s_i(x)$ polynomials
in $x$ of degree at most $d$. $h(x)$ is a multivalued function with at
most $d$ branches, and it may have singularities - poles and ramification
points with finite and infinite values. These singularities appear at zeroes
of $s_d(x)$ and at zeroes of the discriminant $\Delta(x)$ of the equation
(5.1). Consequently, their number does not exceed $r=r(d)=d+d(2d-1)=2d^2$.
Let us denote the singular points of $h(x)$ by $\Sigma=\{x_1,\dots,x_r\}$.

Equivalently, we can consider an algebraic curve $Y\subset {\mathbb C}^2$
defined by the equation $Q(x,y)=0,$ and its projection $\pi_1$ to the
$x$-axis in ${\mathbb C}^2$. The algebraic function $y=h(x)$ is then a
(multivalued) inversion of $\pi_1$.

For a complex polynomial $P(x,y)$ of degree $d_1$ consider $g(x)=P(x,h(x))$.
$g(x)$ is once more an algebraic function. Its singularities are the same
as for $h(x)$. Fix a certain non-singular point $x_0$ and denote $R$ the
distance of $x_0$ to the nearest singularity $x_i$. Now let us fix a regular
branch $\tilde g(x)$ of $g(x)$ at $x_0$. $\tilde g(x)$ can be continued as
a univalued regular function to the disk $D_R(x_0)$ of radius $R$ centered
at $x_0$.

\bl For any $s=0,1,\dots,$ \ $\tilde g(x)$ is a $(s,p)$-valent function in
$D_R(x_0)$, with $p=d \ \max(s,d_1)$.\el\pr For any polynomial $R(x)$ of
degree $s$ the equation $\tilde g(x)=R(x)$ is equivalent to
$\tilde P(x,h(x))=0$, where $\tilde P(x,y)=P(x,y)-R(x)$. The degree of
$\tilde P(x,y)$ is equal to $\kappa=\max(s,d_1)$. Now the equation
$\tilde P(x,h(x))=0$ is equivalent to the system $\tilde P(x,y)=0, \ Q(x,y)=0.$
By Bezout theorem the number of solutions of this system does not exceed
$d\kappa$. $\square$

\smallskip

Now we can state the Remez inequality for $\tilde g$. Let $x_0, R$ be as above.

\bt For $Z\subset D_{R_1}(x_0), \ R_1 < R,$ assume that $c_{d_1}(Z)=c>0$.
Then on each $D_{R'}(x_0)$ with $R'<R$ the function $\tilde g$ satisfies
\be \max_{D_{R'}(x_0)} \vert \tilde g(x) \vert \leq \sigma(R,R_1,R')^{dd_1}
({6eR \over c})^{d_1} \max_{Z}\vert \tilde g(x) \vert,\ee where
$\sigma(R,R_1,R')= \sigma({R_1\over R},{R'\over R}),$ in notations of
Theorem 4.1.\et\pr The number of zeroes of $\tilde g$ in $D_R(x_0)$ is at most
$d_1$. Indeed, substituting $y=0$ into $P(x,y)=0$ we get a polynomial
equation $P(x,0)=0$ in one variable, of degree at most $d_1$, which has at most
$d_1$ solutions. So taking $s=d_1$ we conclude by Lemma 5.1 that
$\tilde g(x)$ is a $(d_1,p)$-valent function in
$D_R(x_0)$, with $p=d \ \max(d_1,d_1)=dd_1.$ Now the required result follows
directly from Theorem 4.1, as we rescale the disk $D_R(x_0)$ to $D_1$. $\square$.

\medskip

Based on Theorem 5.1 one can easily produce a ``global" Remez inequality for
polynomials on algebraic curves. Let a closed set $Z\subset {\mathbb C}\setminus \Sigma$
be given. We assume that $Z\subset D$ for a certain open disk
$D\subset {\mathbb C}\setminus \Sigma$, which we fix. Consider now all the branches
of the algebraic function $y=h(x)$ as above. We fix a certain (regular and
univalued) branch $\hat h(x)$ of $h$ over $D$. Now let us fix a point
$x_0\in {\mathbb C}\setminus \Sigma$ and one of the branches $\bar h(x)$ of
$h$ over a neighborhood of $x_0$.

\medskip

The definition of the constant $K(S)$ below depends only on the geometry of
$\Sigma$, on the position of $x_0$ and $Z$ in ${\mathbb C}\setminus \Sigma$, 
and on the ``combinatorial" data: the degrees $d, \ d_1$, the monodromy group 
$G$ of $h$, and on the choice of the branches $\hat h, \bar h$ of the algebraic 
function $h$. Accordingly, we
shall consider ``configurations" $S=\{d,d_1,\Sigma, Z, x_0, G, \hat h, \bar h \}$
and notice that all the constructions below depend only on these configurations.

\medskip

Let us now consider all the finite chains $CH$ (of a variable length $m$)
of open disks $D_j\subset {\mathbb C}\setminus \Sigma, \ j=0,\dots,m$,
and for each $D_j$ a regular branch $h_j(x)$ of $h$ over $D_j$,
satisfying the following conditions:

\smallskip

1. $D_0=D$, and so $Z\subset D_0$ and $h_0(x) = \hat h(x)$ for $x \in D_0$.

\smallskip

2. $x_0\in D_m$ and $h_m(x) = \bar h(x)$ for $x\in D_m$.

\smallskip

3. For each $j=0,\dots,m-1$ the intersection $D_j \cap D_{j+1}$ is not
empty, and the branches $h_j(x)$ and $h_{j+1}(x)$ coincide over
$D_j \cap D_{j+1}$.

\smallskip

For each such chain $CH$ we define the constant $K(CH)$ as follows: fix the
centers of the disks in $CH$ and denote by $R_j$ the radius of the disk
$D_j, \ j=0,\dots,m$.
First we chose for each $j=0,\dots,m$ the radii $R_{j,1}<R_j$ and $R'_j< R_j$
in such a way that the intersection
$D'_j \cap D_{j+1}$ still is not empty, where $D'_j$ is the disk
concentric with $D_j$ of radius $R'_j$. We also require that $x_0\in D'_m$ and
$Z\subset D^1_0$, where $D^1_j$ is the disk
concentric with $D_j$ of radius $R_{j,1}.$ Put $c_j=c_{d_1}(D_{j,1}\cap D'_{j-1})$.

Now we define $K(CH, d,d_1)$ as
\be K(CH,d,d_1)=\min \prod_{j=0}^m
\sigma(R_j,R_{j1},R'_j)^{dd_1}({{6eR_j}\over c_j})^{d_1},\ee
where the minimum is taken over all the choices of the radii 
$R_{j,1}<R_j, R'_j< R_j, \ j=0,1,\dots,m,$
satisfying the conditions above. Finally we give the following definition:

\bd For a configuration
$S=\{d,d_1,\Sigma, Z, x_0, G, \hat h, \bar h \}$ we define the constant $K(S)$ as
$K(S)=\min K(CH,d,d_1),$ where the minimum is taken over all the choices of the
chains $CH$ as above.\ed

\bt Let an algebraic function $h(x)$ be defined by (5.1), with the configuration
$S=\{d,d_1,\Sigma, Z, x_0, G, \hat h, \bar h \}$. Put $c=c_{d_1}(Z)$.
Then for each polynomial $P(x,y)$ of degree $d_1$ and for $g(x)= P(x,h(x))$
the branches $\hat g$ and $\bar g$ corresponding to the branches $\hat h, \bar h$
of $h$ satisfy
\be \vert \bar g(x_0) \vert \leq K(S)
({{6eR_0}\over c})^{d_1} \max_{Z}\vert \hat g(x) \vert.\ee\et\pr For a fixed
chain $CH$ and the radii $R_{j,1},R'_j$ we apply the inequality of Theorem 5.1
subsequently to the couples $D_j,D_{j,1}\cap D'_{j-1}$ for $j=1,\dots,m$,
and then take minimum with respect to all the choices involved. $\square$

The form of the inequality in Theorem 5.2 separates the roles of the
invariant $c(Z)$, the size of the disk $D_0$ in ${\mathbb C}\setminus \Sigma$
containing $Z$, and of the ``global position" of $Z$ and $x$ with respect to
$\Sigma$ and the chosen branches of $h$. It would be important to give more
explicit estimates for the constant $K(S)$. In particular, one apparent
possible improvement relates to the fact that the {\it total number of zeroes
of $h$ on all its branches is bounded by $d_1$.} Hence we can try to ``separate"
zeroes of $h$ between the disks $D_j$. However, the detailed analysis of the
construction above looks rather tricky, and this problem is beyond
the scope of the present paper.

\bibliographystyle{amsplain}

\end{document}